\documentclass[11pt]{amsart}
\usepackage{amsmath}
\usepackage{graphicx}
\usepackage{comment}
\usepackage{enumerate}
\def\~{\tilde}

\def\a{\alpha}
\def\s{\sigma}

\def\b{{\bar b}}

\def\g{{\bar g}}

\def\L{{\cal L}_p}

\def\S{{\cal S}}

\def\A{{\mathcal A}}

\def\S{{\mathcal S}}
\def\I{{\mathcal I}}
\def\centre{{\mathcal Z}}

\def\dim{{\mathrm {dim}}}

\def\ker{{\mathrm {ker}}}
\def\max{{\mathrm {max}}}

\def\charac{\mathrm {char}}

\def\ann{\mathrm {Ann}}
\def\Sym{\frak{S}}

\usepackage{stmaryrd}

\def\l{{\lambda}}
\def\L{{\Lambda}}

\def\g{{\gamma}}
\def\s{{\sigma}}
\def\a{{\alpha}}
\def\b{{\beta}}

\def\Proof{\noindent{\sl Proof.}\ }
\def\qed{{\hfill $\Box$ \medbreak}}
\usepackage{amssymb}
\usepackage{stmaryrd}

\newtheorem{defi}{Definition}[section]
\newtheorem{thm}[defi]{Theorem}
\newtheorem{lem}[defi]{Lemma}
\newtheorem{cor}[defi]{Corollary}

\newtheorem{eg}[defi]{Example}

\newtheorem{prop}[defi]{Proposition}

\newtheorem*{probA}{Problem A}
\newtheorem*{probB}{Problem B}

\newtheorem*{probC}{Problem C}

\begin{document}

\title[Polynomials that are not quite an identity]{On polynomials that are not quite an identity on an associative algebra}

\author{Eric Jespers, David Riley and Mayada Shahada}
\thanks{The authors acknowledge support from  Onderzoeksraad of Vrije Universiteit, Fonds
voor Wetenschappelijk Onderzoek (Vlaanderen) and NSERC of Canada.}

\begin{abstract} 
Let $f$ be a polynomial in the free algebra over a field $K$, and let $A$ be a $K$-algebra.   
We denote by $\S_A(f)$, $\A_A(f)$ and $\I_A(f)$, respectively, the `verbal' subspace, subalgebra, and ideal, in $A$, 
generated by the set of all $f$-values in $A$.  We begin by studying the following problem: if $\S_A(f)$ is finite-dimensional, is it true that $\A_A(f)$ and $\I_A(f)$ are also finite-dimensional? We then consider the dual to this problem for `marginal' subspaces that are finite-codimensional in $A$. If $f$ is multilinear, the marginal subspace, $\widehat{\S}_A(f)$, of $f$ in $A$ is the set of all elements $z$ in $A$ such that $f$ evaluates to 0 whenever any of the indeterminates in $f$ is evaluated to $z$. We conclude by discussing the relationship between the finite-dimensionality of $\S_A(f)$ and the finite-codimensionality of $\widehat{\S}_A(f)$.
\end{abstract}

\maketitle

\section{Verbal subspaces, subalgebras and ideals}

Throughout this paper, the term `algebra' will be reserved for a not necessarily unital associative algebra $A$ over a fixed base field $K$ of characteristic $p\ge0$.  We shall use $A_1$ to indicate its unital hull.  

\begin{defi}
Let $A$ be an algebra, and let $f=f(x_1,\ldots,x_n)$ be a polynomial in the free algebra $K\langle X\rangle$ on the set $X= \{ x_1,x_2,\ldots \}$.  
We shall denote by $\S_A(f)$, $\A_A(f)$ and $\I_A(f)$, respectively, the subspace, subalgebra and ideal in $A$ 
generated by the set of all {\em $f$-values} in $A$:
\begin{center}
$\{f(a_1,\ldots,a_n):\,a_1,\ldots,a_n \in A\}.$
\end{center}  
We shall call the subspace $\S_A(f)$ the {\em verbal subspace of $A$ generated by $f$}, and so forth.
\end{defi}

It is not hard to construct examples showing that the obvious inclusions
\begin{center}
 $\S_A(f)\subseteq \A_A(f)\subseteq \I_A(f)$
\end{center}
are sometimes strict:
\begin{eg}
\label{first example}
\begin{enumerate}
\item Let $A=M_2(K)$, the algebra of all $2 \times 2$ matrices over $K$, and let $f=[x_1,x_2]=x_1x_2-x_2x_1$.
Then $\S_A(f)=sl_2(K)$, the set of all $2 \times 2$ matrices with trace zero, while $\A_A(f)=\I_A(f)=M_2(K)$.
\item Let $E=E_0 \oplus E_1$ denote the Grassmann algebra of an infinite-dimensional vector space over a field of characteristic $p \neq 2$, and let $f=[x_1,x_2]$. Then $\S_E(f)=\A_E(f) \subseteq E_0$, and yet $\I_E(f) \nsubseteq E_0$.
\item Let $A$ be the polynomial algebra $K[t]$ over a
perfect field $K$ of characteristic $p > 0$, and let $f=x_1^p$. Then $t^{p+1} \in \I_A(f)$, whereas $t^{p+1} \notin \S_A(f)=\A_A(f)=K[t^p]$.
\end{enumerate}
\end{eg}

Even though we have just seen that corresponding verbal subspaces, subalgebras and ideals can be distinct, it is still not clear how different they can be.  Along this vein, we intend to study the following natural problem: 

\begin{probA} Let $A$ be an algebra, and let $f$ be a polynomial in $K\langle X\rangle$.
\begin{enumerate}
\item If $\S_A(f)$ is finite-dimensional, then under what conditions on $A$ and $f$ is $\A_A(f)$ necessarily finite-dimensional? 
\item If $\S_A(f)$ is finite-dimensional,  then under what conditions on $A$ and $f$ is $\I_A(f)$ necessarily finite-dimensional? 
\item If $\A_A(f)$ is finite-dimensional,  then under what conditions on $A$ and $f$ is $\I_A(f)$ necessarily finite-dimensional?
\end{enumerate}
\end{probA}

\begin{eg} \label{example1} Let $p$ be any prime number, let $f=x_1^p$, and let $KG$ be the group algebra of an infinite elementary Abelian $p$-group $G$ over a field $K$ of characteristic $p$. Then
$\S_{KG}(f)= \A_{KG}(f)=K$
is 1-dimensional while $\I_{KG}(f)=KG$ is infinite-dimensional. 
\end{eg} 

In light of the preceding example, it is reasonable to restrict our attention in parts (2) and (3) in Problem A to polynomials $f$ that are linear in at least one indeterminate in the case when the characteristic $p$ is positive.
Before revising our formulation of Problem A, we introduce some additional notation in order to simplify the exposition.

\begin{defi}
Let $f$ be any polynomial in the free algebra $K\langle X\rangle$. 
\begin{enumerate}
\item We shall say that $f$ is {\em algebraically concise} if, for all algebras $A$, $\A_A(f)$ is finite-dimensional whenever $\S_A(f)$ is finite-dimensional.
\item We shall say that $f$ is {\em ideally concise} if, for all algebras $A$, $\I_A(f)$ is finite-dimensional whenever $\S_A(f)$ is finite-dimensional.\end{enumerate}
\end{defi}

It seems plausible to the authors that the following two problems implicit within Problem A have positive solutions.  These questions form the primary focus of our first area of inquiry.  By a slight abuse of notation, we shall refer to them in the sequel by their original labels.

\begin{enumerate}
\item[A(1):] {\em Is every polynomial $f$ in $K\langle X\rangle$ algebraically concise?}
\item[A(2):] {\em Is every polynomial $f$ in $K\langle X\rangle$ ideally concise provided either $p=0$ or $p>0$ and $f$ in linear in at least one indeterminate?}
\end{enumerate} 

We claim first that, unlike the case when the characteristic is positive, restricting our attention to multilinear polynomials $f$ in the case when the characteristic is zero is not actually a restriction at all. Indeed, this assertion is direct consequence of the following proposition.

\begin{prop} \label{multilinear} Let $K$ be a field of characteristic zero. Then a polynomial $f$ in $K\langle X\rangle$ is algebraically concise (respectively, ideally concise) if and only if each of its multilinear consequences is algebraically concise (respectively, ideally concise).
\end{prop} 

Before proving Proposition \ref{multilinear}, we first reduce to the case of homogeneous polynomials. Recall that a polynomial $f$ is called {\em homogeneous} whenever the degree of each indeterminate appearing in $f$ does not vary among the monomials in $f$.

\begin{prop} \label{homogeneous} Let $K$ be a field of arbitrary characteristic, and suppose that $f$ is any polynomial  in $K\langle X\rangle$ with the property that $|K| \geq m+1$, where $m$ is the maximum degree of any indeterminate appearing in $f$.   Then a  polynomial $f$ in $K\langle X\rangle$ is algebraically concise (respectively, ideally concise) if and only if each of its homogenous components is algebraically concise (respectively, ideally concise).
\end{prop}

The proofs of Propositions \ref{multilinear} and \ref{homogeneous} will be carried out in Section 2. 

Our most complete result is a partial solution to Problem A(1):

\begin{thm}
\label{main} Let $f$ be any polynomial in the free algebra over a field $K$, and let $m$ be the maximum degree of an indeterminate appearing in $f$.
\begin{enumerate}
\item If $|K|\ge m+1$, then $f$ is algebraically concise.
\item If $|K|\ge m$ and $f$ is homogeneous, then $f$ is algebraically concise.
\end{enumerate}
In particular, if either $K$ is infinite or $K$ is finite and $f$ is multilinear, then $f$ is algebraically concise.
\end{thm}

The following is a direct consequence of Theorem \ref{main} and Example \ref{example1}.

\begin{cor} Let $K$ be any field of positive characteristic $p$. Then $f=x_1^p$ is algebraically concise but not ideally concise.
\end{cor}

The following interesting special cases of Problem A(1) remain open for (small) finite base fields.
{\em
\begin{itemize} 
\item Are all polynomials of the form $x_1^{p^n}$ algebraically concise?
\item Are all Engel polynomials $[x_1,\underbrace{x_2,\ldots,x_2}_n]$ algebraically concise? 
\end{itemize}
}

 We shall refer to (left-justified) multilinear Lie monomials of the form $[x_1,x_2,\ldots,x_{n}]$ as being {\em simple} of {\em length} $n\ge1$. Multilinear Lie monomials formed by arbitrary bracketing are called {\em outer Lie commutators}.
The sequence of {\em derived series commutators} are defined recursively by $g_1=[x_1,x_2]$
and, for each $n\ge2$,
\begin{center}
$g_n=[g_{n-1}(x_1,\ldots,x_{2^{n-1}}),g_{n-1}(x_{2^{n-1}+1},\ldots,x_{2^n})].$
\end{center}

Our proof of Theorem \ref{main} depends on the following partial solution to Problem A(2):
\begin{thm}
\label{thmB}
Let $f$ be any polynomial in the free algebra $K\langle X\rangle$.  Then the polynomial $g$ defined by
$$g(x_1,\ldots,x_{n+1})=[f(x_1,\ldots,x_n),x_{n+1}]$$
is ideally concise.
\end{thm}
 
The following examples are noteworthy special cases of Theorem \ref{thmB}.

\begin{cor} Let $K$ be any base field, and let $n$ be a positive integer. Then the following statements hold.
\begin{enumerate}
\item The simple Lie commutator $[x_1,x_2,\ldots,x_{n+1}]$ is ideally concise.
\item The outer Lie commutator $[g_n,x_{2^n+1}]$ is ideally concise.  In particular, the Lie centre-by-metabelian polynomial $[[[x_1,x_2],[x_3,x_4]],x_5]$ is ideally concise.\end{enumerate}
\end{cor}

Recall that the {\em standard polynomial} of degree $n$ is given by
\begin{center} $s_{n}(x_1,\ldots,x_{n})=\sum_{\s\in \Sym_{n}}(-1)^{\s}x_{\s(1)}\cdots x_{\s(n)},$\end{center} 
where $\Sym_{n}$ denotes the symmetric group of degree $n$.
We remark that Theorem \ref{thmB} leaves open the following interesting special cases of Problem A(2). {\em
\begin{itemize}
\item Are all the standard polynomials $s_n$ ideally concise?
\item Are the derived series commutators, starting with the Lie metabelian commutator $g_2=[[x_1,x_2],[x_3,x_4]]$, ideally concise?
Are all outer Lie commutators ideally concise?
\item Are all Engel polynomials ideally concise?
\end{itemize}
}

The proofs of Theorems \ref{main} and \ref{thmB} will be carried out in Section 2. We will also provide there a complete positive solution to our original Problem A(2) in the special case when $A$ is a group algebra. 

We close this section by reframing Problem A(2) within the context of polynomial identity theory.

\begin{defi}
Let $A$ be an algebra, and let $f(x_1,\ldots,x_n)$ be any polynomial in $K\langle X\rangle$. If $f(a_1,\ldots,a_n)=0$, for all $a_1,\ldots,a_n$ in $A$,
then $f$ is called a {\em polynomial identity} on $A$. Polynomial identities are called PIs for short.
 It is customary to denote 
by $T(A)$ the set of all PIs on $A$. 
If $T(A)$ is nontrivial, then $A$ is called a {\em PI-algebra}. 
If $\S_A(f)$ is finite-dimensional, we shall say that $f$ is {\em almost a polynomial identity on} $A$;   
we shall denote by $\widetilde{T}(A)$ the set of all such polynomials $f$.  Clearly, $T(A)\subseteq \widetilde{T}(A)$.
\end{defi}

As a consequence of the following lemma, we may assume that every algebra $A$ in Problem A is an infinite-dimensional PI-algebra.
 
\begin{lem}
\label{PI} Let $A$ be an algebra. If there exists a nonzero polynomial $f$ such that $f$ is almost a PI for $A$, then $A$ is a PI-algebra. In other words,  if $\widetilde{T}(A)\neq0$ then $T(A)\neq0$. 
\end{lem}

\Proof Suppose that $\dim(\S_A(f))=t<\infty$. Then the antisymmetry of the standard polynomial $s_{t+1}$,
together with the linear dependence of any $(t+1)$-many $f$-values in $A$, forces $A$ to satisfy the nontrivial multilinear polynomial identity $g$, where $g$ is defined by
\begin{center}
$g=s_{t+1}(f(x_1,\ldots, x_n),\ldots, f(z_1,\ldots, z_n))$.
\end{center}
\qed

A subspace of the free algebra is called a {\em $T$-subspace} whenever it is closed under the algebra endomorphisms of the free algebra.  The notion of $T$-subspace extends to arbitrary algebras in the obvious way. The terminology is also applied to subalgebras and ideals; for example, $T(A)$ is clearly a $T$-ideal of the free algebra, for every algebra $A$. The proof of our next observation is straightforward and thus omitted. 

\begin{lem}
\label{Tspaces} Let $A$ be an algebra.  
Then the following statements hold.
\begin{enumerate}
\item $\S_A(f)$, $\A_A(f)$ and $\I_A(f)$ are $T$-subspaces of $A$, for each $f\in K \langle X\rangle$.
\item $\widetilde{T}(A)$ is a $T$-subalgebra of $K \langle X\rangle$. 
\item $\widetilde{T}(A)$ is a $T$-ideal of $K \langle X\rangle$ if and only if, for every $f\in \widetilde{T}(A)$, $\I_A(f)$ is finite-dimensional. 
\end{enumerate}
\end{lem}

As a consequence of part (3) of Lemma \ref{Tspaces}, the characteristic zero case of Problem A(2) is precisely equivalent to:
\begin{itemize}
\item Is $\widetilde{T}(A)$ a $T$-ideal of $K \langle X\rangle$, for every algebra $A$?
\end{itemize}

In Section 3, we shall introduce marginal subspaces, which are a kind of dual notion to verbal subspaces, and then study the relationship between marginal subspaces, subalgebras and ideals.  Finally, in Section 4, we will discuss the relationship between verbal and marginal spaces. 

All the problems that we propose to study, together with the terminology we use, were inspired by Philip Hall's seminal work on verbal and marginal subgroups beginning in 1940 (see \cite{PH} for his collected works). These problems stimulated a broad stream of deep and interesting research in group theory that remains highly active today. In particular, Hall conjectured that every word in the free group is concise. Although many partial positive solutions of this conjecture have been since proved, the conjecture in its most general form was eventually refuted by Ivanov in 1989 (see \cite{Iva}).  

Ian Stewart made a study of the natural Lie-theoretic analogue of Hall's conjecture in \cite{STW}.  Given the relative complexity of the corresponding problems in the categories of groups and associative algebras, Stewart proved the following surprisingly strong result (Corollary 3.2 in \cite{STW}): if $K$ is an infinite field, $L$ is a Lie algebra over $K$, and $f$ is any polynomial in the free Lie algebra over $K$, then $\S_L(f)=\A_L(f)=\I_L(f)$ is a characteristic ideal in $L$.  As a trivial consequence of this remarkable fact,  in the category of Lie algebras (over an infinite field), every polynomial is ideally concise.

\section{Proof of Section 1 results}

The same argument used to prove Corollary 2.2 in  \cite{V} when applied to (associative) algebras over `sufficiently large' fields yields the following lemma.

\begin{lem}
\label{lem1}
Let $K$ be a field, and suppose that $f$ is any polynomial  in $K\langle X\rangle$ with the property that $|K| \geq m+1$, where $m$ is the maximum degree of any indeterminate appearing in $f$. If $f=\sum_{i=1}^r f_i$ is the homogeneous decomposition of $f$, then, for every algebra $A$, we have 
$$\S_A(f)=\S_A(f_1)+\cdots+\S_A(f_r).$$ 
\end{lem}

Recall that every algebra $A$ can be viewed as a Lie algebra via its Lie bracket $[a,b]=ab-ba$. The following lemma is well-known and straightforward to verify.

\begin{lem}
\label{id}
The following identities hold for all algebras $A$.
\begin{enumerate}
\item Adjoint maps are derivations; in other words, for all $a,b,c \in A$,
$$[ab,c]=a[b,c]+[a,c]b.$$
\item The semi-Jacobi identity holds; namely, for all $a,b,c \in A$,
$$[ab,c]=[a,bc]+[b,ca].$$
\end{enumerate}
\end{lem}

The following result will play a key role throughout this section.

\begin{thm}
\label{Lie}
Let $f$ be any polynomial in the free algebra $K\langle X\rangle$, let $m$ be the maximum degree of any indeterminate appearing in $f$, and let $A$ be an algebra. If either $|K| \ge m+1$ or $|K|\ge m$ and $f$ is homogeneous, then $\S_A(f)$ is closed under derivations, so that $\S_A(f)$ is a Lie ideal of $A$. 
In particular, if either $K$ is infinite or $f$ is multilinear, then $\S_A(f)$ is a Lie ideal of $A$.
\end{thm}

\Proof 
If $K$ is infinite, it follows from Theorem 3.1 in \cite{STW} that $\S_A(f)$ is closed under derivations of $A$.  Since part (1) of Lemma \ref{id} says that adjoint maps are associative derivations, it follows that $\S_A(f)$ is a Lie ideal in $A$.  In the case when $K$ is finite and $|K| \ge m+1$, applying Lemma \ref{lem1} allows us to assume that $f$ is homogeneous. The result now follows just as before, but this time using Lemma 3.4 in \cite{STW} in place of Theorem 3.1 in \cite{STW}.
\qed

\begin{lem}
\label{reduction}
Let $J=W_1+\cdots+W_t$ be any subspace decomposition of a Lie ideal $J$ of an (associative) algebra $A$.  Then the following statements hold:
\begin{enumerate}
\item The unital subalgebra of $A_1$ generated by $J$ is precisely $(S_1)_1\cdots (S_t)_1$, where each $S_i$ is the subalgebra of $A$ generated by $W_i$.
\item The (two-sided) ideal in $A$ generated by $J$ is $W_1A_1+\cdots+W_tA_1$.
\end{enumerate}
\end{lem}
\Proof
First notice that, if $w_i\in W_i$ and $w_j\in W_j$, then $w_jw_i=w_iw_j+[w_j,w_i]\in w_iw_j+J=w_iw_j+W_1+\cdots + W_t$.  Thus, longer products of homogeneous elements can be reordered modulo shorter products. Part (1) now follows by an induction on the length of the products.  To prove Part (2), we need only observe that $aw_i=w_ia+[a,w_i]\in w_ia+J$, for all $a\in A$ and $w_i\in W_i$.
\qed

In order to deduce Proposition \ref{homogeneous} from our preliminary results, suppose that $K$ is sufficiently large (as described in its hypothesis) for a given (non-homogeneous) polynomial $f$ with homogeneous decomposition $f=\sum_{i=1}^r f_i$.  Then, by Lemma \ref{lem1} and Theorem \ref{Lie}, $\S_A(f)$ is a Lie ideal of $A$ that coincides precisely with the vector space sum of the Lie ideals $\S_A(f_1),\ldots,\S_A(f_r)$. Consequently,  Lemma \ref{reduction} implies the following:

\begin{prop} \label{homogeneousdecomp}
Let $K$ be a field, and suppose that $f$ is any polynomial  in $K\langle X\rangle$ with the property that $|K| \geq m+1$, where $m$ is the maximum degree of any indeterminate appearing in $f$. If $f=\sum_{i=1}^r f_i$ is the homogeneous decomposition of $f$, then, for every algebra $A$, we have 
\begin{enumerate}
\item $\S_A(f)=\S_A(f_1)+\cdots +\S_A(f_r)$; 
\item $(\A_A(f))_1= (\A_A(f_1))_1\cdots (\A_A(f_r))_1$; and,
\item $\I_A(f)=\I_A(f_1)+\cdots +\I_A(f_r)$.
\end{enumerate}
\end{prop}

Proposition \ref{homogeneous} is a straightforward consequence of Proposition \ref{homogeneousdecomp}.
 
Proposition \ref{homogeneous} allows us to assume that $f$ is homogeneous in Proposition \ref{multilinear}.  Applying the process of multilinearization to the homogeneous polynomial $f$ shows that
$$\S_A(f) \supseteq\sum_{j=1}^s \S_A(g_j),$$
where $\{g_1,\ldots,g_s\}$ is a finite set of multilinear consequences of $f$ with the same total degree as $f$. When the characteristic is zero, a simple Vandermonde argument shows that this inclusion can be taken to be an equality.  Theorem \ref{Lie} and Lemma \ref{reduction} now proves Proposition \ref{multilinear}.


The following result is a constructive version of Theorem \ref{thmB}.

\begin{thm}
\label{thm1}
Let $f$ be any polynomial in $K\langle X\rangle$, let $m$ be the maximum degree of any indeterminate appearing in $f$, and suppose that $A$ is an algebra such that $\dim(\S_A(g)) =t$ is finite, where the polynomial $g$ is defined by
\begin{center}
$g(x_1,\ldots,x_{n+1})=[f(x_1,\ldots,x_n),x_{n+1}].$
\end{center}
Then
\begin{center}
$\dim(\I_A(g)) \leq t(t^2+1)^2.$
\end{center}
In the case when either $|K| \geq m+1$ or $|K|\ge m$ and $f$ is homogeneous, this bound can be sharpened to
\begin{center}
$\dim(\I_A(g)) \leq t(t^2+1)$.
\end{center} 
In particular, the sharper bound holds when $K$ is infinite or $f$ is multilinear.
\end{thm}

\Proof We can choose a basis $\{\bar{g}_1,\ldots, \bar{g}_t\}$ of $\S_A(g)$ such that each basis element $\bar{g}_i$ has the form $\bar{g}_i=[\bar{f}_i,a_i]$ for some evaluation $\bar{f}_i$ of $f$ in $A$. Denote by $C_A(\bar{f}_i)$ the centralizer of $\bar{f}_i$ in $A$.  Then, for each $i$, $\dim (A/C_A(\bar{f}_i)) \leq t$ (consider the kernel of the linear map $A \rightarrow \S_A(g): a \mapsto [\bar{f}_i,a] $). Hence, $\dim (A/C) \leq t^2$, where 
\begin{center}
$C= \bigcap_{i=1}^t C_A(\bar{f}_i)$.
\end{center}
Thus, there exist $b_1,b_2,\ldots,b_{t^2} \in A$ such that $A_1=K1+Kb_1+\cdots+Kb_{t^2} + C$.

We claim that $\S_A(g)C \subseteq \S_A(g)$. Indeed, by part (1) of Lemma \ref{id}, for each $i=1,\ldots,t$ and $c \in C$, we have
\begin{align}
\nonumber [\bar{f}_i,a_i]c & = [\bar{f}_i,a_ic]-a_i[\bar{f}_i,c]\\
\nonumber &=[\bar{f}_i,a_ic] \in \S_A(g).
\end{align}
Similarly, $C\S_A(g) \subseteq \S_A(g)$.  Therefore, 
\begin{align}
\nonumber \I_A(g) &= A_1\S_A(g)A_1\\
\nonumber &= (K1+Kb_1+\cdots+Kb_{t^2})\S_A(g)(K1+Kb_1+\cdots+Kb_{t^2}),
\end{align}
so that $\dim(\I_A(g)) \leq t(t^2+1)^2$. 

We can improve this bound somewhat in the case when $|K| \geq m+1$ or $|K|\ge m$ and $f$ is homogeneous. Indeed, in this case, Theorem \ref{Lie} informs us that $\S_A(f)$ is a Lie ideal of $A$.  By the same reasoning, $\S_A(g)$ is also a Lie ideal. Consequently, we have 
\begin{align}
\nonumber \I_A(g) &= \S_A(g)A_1\\
\nonumber &= \S_A(g)(K1+Kb_1+\cdots+Kb_{t^2}),
\end{align}
so that $\dim(\I_A(g)) \leq t(t^2+1)$.
\qed
 
 Theorem \ref{main} is a direct consequence of the following constructive result.
 
\begin{thm}
\label{main1}
Let $f$ be any polynomial in $K\langle X\rangle$ with homogeneous components $\{f_\l\}_{\l  \in \Lambda}$, and let $r=\max\{r_\l\,|\,\l \in \Lambda\}$, where $r_\l$ is the minimum degree of any indeterminate appearing in $f_\l$.  Suppose that either $|K| \geq m+1$ or $|K|\ge m$ and $f$ is homogeneous, where $m$ is the maximum degree of any indeterminate appearing in $f$, and suppose that $A$ is an algebra such that $\dim(\S_A(f)) = t$ is finite.
 Then
$$ \dim(\A_A(f)) \leq (r^t-1) + tr^{t-1} +t(t^2+1).$$
In particular, in the case when $f$ is multilinear, the following bound holds:
$$ \dim(\A_A(f))\leq t+t(t^2 + 1).$$
\end{thm}
\Proof 
Suppose that either $|K| \geq m+1$ or $|K|\ge m$ and $f$ is homogeneous, and set $g(x_1,\ldots,x_{n+1})=[f(x_1,\ldots,x_n),x_{n+1}]$.
Then, according to Lemma \ref{lem1} and Theorem \ref{Lie}, 
$$\S_A(f)=\sum_{\l \in \Lambda} \S_A(f_\l)$$
is a Lie ideal of $A$.  
Let $z_1,\ldots,z_t$ be a basis of $\S_A(f)$ chosen with the property that each $z_i$ lies in $\S_A(f_\l)$, for some $\l \in \Lambda$.  Since $\S_A(f)$ is a Lie ideal of $A$, we have $\S_A(g) \subseteq \S_A(f)$, so that $\dim(\S_A(g)) \leq t$. Thus, by Theorem \ref{thm1}, we may replace $A$ by $A/\I_A(g)$ to assume that $g$ is a polynomial identity satisfied by $A$; in other words, we may assume that $\S_A(f)$ is contained in the centre of $A$. It follows that 
\begin{center}
$\S_A(f_\l)z^{r_\l} \subseteq \S_A(f_\l) \subseteq \S_A(f)$,
\end{center} 
for each $z \in \A_A(f)$ and $\l \in \Lambda$.

We claim that $$\A_A(f)=\sum K z_1^{\b_1} \cdots z_t^{\b_t},$$
where each $0\le \b_i\le r$, not every $\b_i$ is zero, and at most one $\b_i=r$. Clearly the righthand subspace is contained in the lefthand subspace.  To prove the reverse inclusion, it suffices to show that every product of the form
$$\pi=z_1^{\a_1} \cdots z_t^{\a_t},$$
where each $\a_i\ge0$ and $\a_1+\cdots+\a_t\ge1$, lies in the righthand side.  

Observe first that, if $\a\ge r$ and $\l\in\L$, then, by the division algorithm, $\a = qr_\l + \b$, where $q\ge0$ and $0 \leq \b \leq r_\l-1\le r-1$.  Thus,
$$\S_A(f_\l)z_i^{\a} = \S_A(f_\l)(z_i^{q})^{r_{\l}}z_i^{\b} \subseteq \S_A(f_\l) z_i^{\b}.$$

Consider now the case when $\pi=z_i^{\a_i}$, where $\a_i\ge r+1$. Let $\l\in\L$ be such that $z_i\in\S_A(f_\l)$. Then, by the division algorithm, $\a_i-1 = q_ir_\l + \b_i$, where $q_i\ge0$ and $0 \leq \b_i \leq r_\l-1\le r-1$.  Therefore, as shown above, 
$$\pi=z_i^{\a_i} = z_iz_i^{\a_i-1} \subseteq \S_A(f_\l)z_i^{\a_i-1} \subseteq\S_A(f_\l) z_i^{\b_i} \subseteq \S_A(f) z_i^{\b_i}.$$
This proves the claim in this case.

We can now assume that at least two exponents $\a_i$ in the monomial $\pi=z_1^{\a_1} \cdots z_t^{\a_t}$ are positive with one at least $r$. By commutativity, we may assume $\a_1\ge r$ and $\a_2\ge1$.  Let $\l\in\L$ be such that $z_2\in\S_A(f_\l)$. Then, as shown above,  
$$z_1^{\a_1}z_2^{\a_2} \in \S_A(f_\l)z_1^{\a_1}z_2^{\a_2-1} \subseteq \S_A(f_\l) z_1^{\b_1}z_2^{\a_2-1} \subseteq \S_A(f) z_1^{\b_1}z_2^{\a_2-1},
$$
for some  $0\le\b_1\le r-1$.  Therefore,  
$$\pi=z_1^{\a_1} \cdots z_t^{\a_t}\in \S_A(f) z_1^{\b_1}z_2^{\a_2-1}z_3^{\a_3}\cdots z_t^{\a_t}.$$
Repeating this sort of argument shows that
$$\pi=z_1^{\a_1} \cdots z_t^{\a_t}\in  \S_A(f) z_1^{\b_1}\cdots z_t^{\b_t},$$
where each $0\le\b_i\le r-1$, as required.

It follows from the claim that
$$ \dim(\A_A(f))\leq (r^t -1) + tr^{t-1} + \dim (\I_A(g)),$$
as required. The rest now follows from Theorem \ref{thm1}.
\qed

The proof of Theorem \ref{main1} also yields the following result that may be interesting in its own right.

\begin{cor} If $A$ is an algebra and $f$ is a multilinear polynomial, then
$$
\A_A(f) \subseteq \S_A(f)+\I_A(g),
$$
where $g(x_1,\ldots,x_{n+1})=[f(x_1,\ldots,x_n),x_{n+1}].$
\end{cor}

Problem A is settled in the case of group algebras by the following result. 

\begin{thm}\label{groupalgebras} Let $KG$ be the group algebra of an infinite group $G$ over a field $K$ of characteristic $p\ge0$, and let $f$ be a polynomial over $K$ with the property that $f$ is almost a PI for $KG$. If either  $p=0$ or $p>0$ and $f$ is linear in at least one indeterminate, then $f$ is an actual polynomial identity for $KG$ (so that  $\S_{KG}(f)=\A_{KG}(f)=\I_{KG}(f)=0$ are finite-dimensional).
\end{thm}

The authors are indebted to Donald Passman for sharing with them the following greatly simplified proof of Theorem \ref{groupalgebras}.
\medbreak

\Proof
We can assume that $f(x_1,\ldots, x_n)$ is linear in the last indeterminate, $x_n$. Let $a_1,\ldots, a_{n-1}$ be arbitrary elements in $KG$. Then $$f(a_1,\ldots, a_{n-1}, x_n)=\a_1x_n\b_1+\cdots+\a_sx_n\b_s,$$
for some fixed elements $\a_1,\ldots,\a_s,\b_1,\ldots,\b_s$ in $KG$.
Since $\S_{KG}(f)$ is finite-dimensional, there are only finitely many group elements in its support. It follows that there are only finitely many elements $g$ in $G$ such that $f(a_1,\ldots, a_{n-1}, g)$ is nontrivial. Let $T$ be the complement of these elements in $G$, so that 
$$f(a_1,\ldots, a_{n-1}, g)=\a_1g\b_1+\cdots+\a_sg\b_s=0,$$
for all $g\in T$. By Lemma 4.2.3(iv) in \cite{DSP}, $T$ is clearly very large in $G$, where \lq very large\rq\ is as defined in the paragraph preceding the statement of the lemma.  Consequently, by Lemma 4.2.4 in \cite{DSP}, $f(a_1,\ldots, a_{n-1}, g)=0$, for all $g$ in $G$.  Since $a_1,\ldots, a_{n-1}$ were arbitrary, this means $f=0$ is actually a polynomial identity for $KG$.
\qed

Because Engel polynomials are linear in one indeterminate, we have the following interesting consequence of Theorem \ref{groupalgebras}.

\begin{cor} Let $KG$ be the group algebra of an infinite group $G$ over an arbitrary field $K$ of characteristic $p\ge0$.  If 
\begin{center}
$[x,\underbrace{y,\ldots,y}_n]=0$
\end{center}
is almost a polynomial identity for $KG$, then it is an actual polynomial identity for $KG$. 
\end{cor}

We remark that group algebras satisfying nontrivial polynomial identities were characterized (in terms of the group structure of $G$) by Isaacs and Passman (see Section 5.3 in \cite{DSP}), while group algebras satisfying an Engel identity were later characterized by Sehgal  (see Section V.6 in \cite{SKS}).

\section{Marginal subspaces, subalgebras and ideals}

In this section, we study the natural dual to a verbal subspace.

\begin{defi}
Let $A$ be any algebra, and let $f(x_1,\ldots,x_n)$ be any polynomial in the free algebra $K \langle X \rangle$. 
\begin{enumerate}
\item We shall call an element $z \in A$ with the property that
\begin{center}
$f(b_1,\ldots,b_{i-1},b_i+\a z,b_{i+1},\ldots,b_n)-f(b_1,\ldots,b_n)=0$,
\end{center}
for all choices of $b_1,\ldots, b_n$ in $A$, $\a \in K$, and $i=1,2,\ldots,n$, an {\em eradicator} of $f$ in $A$.  The set of all eradicators of $f$ in $A$ forms a subspace $\widehat{\S}_A(f)$ of $A$ that we shall call the {\em marginal subspace of $A$ with respect to $f$}. 
\item We shall denote by $\widehat{\A}_A(f)$ and $\widehat{\I}_A(f)$ (respectively) the largest subalgebra and ideal of $A$ contained in $\widehat{\S}_A(f)$, which we shall call the {\em marginal subalgebra} and {\em marginal ideal of $A$ with respect to $f$}. 
\end{enumerate}
\end{defi}

In the case when the polynomial $f$ is multilinear, it is easy to see that, $z$ eradicates $f$ in $A$ if and only if
\begin{center}
$f(b_1,\ldots, b_{i-1},z,b_{i+1},\ldots,b_n)=0$,
\end{center}
for all choices of $b_1,\ldots, b_n$ in $A$ and $i=1,2,\ldots,n$. 

We intend to study certain relationships between our various notions of marginal subspaces.  
The following examples are either well-known or straightforward to verify.  

\begin{eg} Let $A$ be an algebra, and let $n$ be a positive integer.
\label{marginalexamples}
\begin{enumerate}
\item If $A$ is commutative and $f$ is a multilinear polynomial, then clearly $\widehat{\S}_A(f)=\widehat{\A}_A(f)=\widehat{\I}_A(f)$.  
\item If $f=x_1 \cdots x_{n+1}$, then $\widehat{\S}_A(f)=\widehat{\A}_A(f)=\widehat{\I}_A(f)=\ann^n(A)$, where $\ann^0(A)=0$ and, for each $n \geq 0$, $\ann^{n+1}(A)$ is given by
\begin{center}
$\ann^{n+1}(A)/\ann^n(A)=\ann(A/\ann^n(A))$,
\end{center}
the two-sided annihilator of the algebra $A/\ann^n(A)$.
\item If $f=[x_1,x_2]$, then $\widehat{\S}_A(f)=\widehat{\A}_A(f)=\centre(A)$, the centre of $A$, while $\widehat{\I}_A(f)=\centre(A)\cap\ann[A,A]$.  For example, let $E$ be the Grassmann algebra of an infinite-dimensional vector space over a field of characteristic $p\neq 2$.  Then $\widehat{\S}_E(f)=\widehat{\A}_E(f)=\centre(E)=E_0$ and $\widehat{\I}_E(f)=0$. 
\item If $f=[x_1, \ldots, x_{n+1}]$, then $\widehat{\S}_A(f)=\widehat{\A}_A(f)=\centre_n(A)$, where $\centre_n(A)$ denotes the $n^{\text{th}}$ term of the ascending central series of $A$ when viewed as a Lie algebra. 
\item Let $f=s_n$ be the standard polynomial of degree $n\ge2$, and let $E$ be the Grassmann algebra over an infinite-dimensional vector space over a field of characteristic zero.  Then $\widehat{\S}_E(f)=\widehat{\A}_E(f)=E_0$, while $\widehat{\I}_E(f)=0$.
\end{enumerate}
\end{eg}

Notice that, in each example, $\widehat{\S}_A(f)=\widehat{\A}_A(f)$; in other words, $\widehat{\S}_A(f)$ is a subalgebra of $A$. In this section, we shall address the following problem, which is helpful to think of as a kind of dual to Problem A in Section 1. 

\begin{probB} Let $A$ be an algebra, and let $f$ be a polynomial in $K\langle X \rangle$.
\begin{enumerate}
\item Under what conditions on $A$ and $f$ is $\widehat{\S}_A(f)$ a subalgebra of $A$? Is this always true?
\item If $A/\widehat{\S}_A(f)$ is finite-dimensional, under what conditions on $A$ and $f$ is $A/\widehat{\A}_A(f)$ necessarily finite-dimensional? Is this always true?
\item If $A/\widehat{\S}_A(f)$ is finite-dimensional, under what conditions on $A$ and $f$ is $A/\widehat{\I}_A(f)$ necessarily finite-dimensional?
Is this always true? 
\item If $A/\widehat{\A}_A(f)$ is finite-dimensional, under what conditions on $A$ and $f$ is $A/\widehat{\I}_A(f)$ necessarily finite-dimensional?
Is this always true? 
\end{enumerate} 
\end{probB}

Problems B(2) and B(3) are, in fact, equivalent since Problem B(4) has a positive solution (in full generality) by the following result of Lee and Liu:

\begin{thm} {\em (\cite{LL})}
\label{LLthm}
Let $A$ be an algebra over a field $K$, and let $B$ be a subalgebra of $A$ such that $\dim(A/B) < \infty$. Then there exists an ideal $I$ of $A$ contained in $B$ such that $\dim(A/I) < \infty$.
\end{thm}

Consequently, if $\widehat{\A}_A(f)$ is finite-codimensional in $A$, then so is $\widehat{\I}_A(f)$.  This means that any positive solution to Problem B(1) is also a positive solution to the other parts of Problem B.

\begin{defi} Let $f$ be a polynomial in the free algebra over $K$. 
\begin{enumerate} 
\item We shall say that $f$ is {\em marginally a PI} on an algebra $A$ whenever $A/\widehat{\S}_A(f)$ is finite-dimensional.
\item We shall call $f$ {\em marginally concise} if, for all algebras $A$, $A/\widehat{\I}_A(f)$ is finite-dimensional whenever $A/\widehat{\S}_A(f)$ is finite-dimensional.
\item We shall call $f$ {\em marginally perfect} if, for all algebras $A$, $\widehat{\S}_A(f)$ is a subalgebra of $A$.
\end{enumerate}
\end{defi}

Observe that $f$ is a PI on $A$ precisely when $\widehat{\S}_A(f)=A$. Also notice that, due to Theorem \ref{LLthm}, there is no need to distinguish between 
being \lq marginally algebraically concise\rq\ and 
\lq marginally ideally concise\rq.  

We intend to demonstrate below that certain classes of polynomials are marginally concise by showing that they are marginally perfect.   
Along these lines, the following corollary of a result of Stewart is worthy of mention in its own right.

\begin{prop} 
Let $A$ be an algebra, let $f$ be a polynomial in the free algebra over $K$, and suppose that either $\charac(K)=0$ or $f$ is multilinear. Then $\widehat{\S}_A(f)$ is invariant under all derivations of $A$; in particular, $\widehat{\S}_A(f)$ is a Lie ideal of $A$.
\end{prop}

\Proof According to Proposition 5.1 in \cite{STW},
$\widehat{\S}_A(f)$ is invariant under all derivations of $A$.  Part (1) of Lemma \ref{id} says that adjoint maps are (associative) derivations. \qed

\begin{lem}
\label{commutator}
The following statements hold in the free $K$-algebra on the set $X=\{w,x_1,x_2,\ldots, x_n \}$. 
\begin{enumerate}
\item For each integer $1 \leq i \leq n$, 
\begin{align}
\nonumber [x_1,\ldots,x_{i-1},x_i,x_{i+1},\ldots,x_n] &=-[x_i,[x_1,\ldots,x_{i-1}],x_{i+1},\ldots,x_n]\\
\nonumber & \in [x_i, \underbrace{K\langle X \rangle,\ldots,K \langle X \rangle}_{n-1}].
\end{align}
\item $[wx_1,x_2,\ldots,x_n]=[w,x_1x_2,x_3,\ldots,x_n]+[x_1,x_2w,x_3,\ldots,x_n]$.
\end{enumerate}
\end{lem}

\Proof
Part (1) follows easily from the Jacobi identity and induction. Part (2) follows from the semi-Jacobi identity (part (2) of Lemma \ref{id}).  
\qed

Our next result provides the first of three partial solutions to Problem B(1). In its proof, we shall use $\bar{f}$ to denote a generic evaluation of a polynomial $f$. Recall that we refer to a left-justified Lie commutator of the form $[x_1,x_2,\ldots,x_{n}]$ as being simple of length $n\ge1$.

\begin{thm}
\label{polyconcise}
Suppose that $f$ is a multilinear polynomial of the form
\begin{center}
$f(x_1,\ldots,x_n) = h_1 \cdots h_m$,
\end{center}
where each polynomial $h_j$ is a simple Lie commutator of length at least one. Then $f$ is marginally perfect.  In particular, simple Lie commutators are marginally perfect.
\end{thm}

\Proof Let $A$ be any algebra, and let $z_1,z_2 \in \widehat{\S}_A(f)$.We need to show that $z_1z_2$ eradicates $f$ in $A$, so that  $\widehat{\S}_A(f)$ is a subalgebra of $A$. Note first that any indeterminate $x_i$ involved in the multilinear polynomial $f$ must fall in exactly one $h_j$. We divide the proof into two cases:\\

\noindent {Case (1):} Suppose that $h_j=x_i$ has length 1. Evaluating $x_i$ by $z_1z_2$ (and the remaining indeterminates with arbitrary elements from $A$) yields:
\begin{align}
\nonumber \bar{f}&=\bar{h}_1 \cdots \bar{h}_{j-1} z_1z_2 \bar{h}_{j+1}\cdots \bar{h}_m\\
\nonumber & = (\bar{h}_1 \cdots \bar{h}_{j-1}z_1\bar{h}_{j+1}z_2\bar{h}_{j+2} \cdots \bar{h}_m)-
(\bar{h}_1 \cdots \bar{h}_{j-1}z_1[\bar{h}_{j+1},z_2]\bar{h}_{j+2} \cdots\bar{h}_m)\\
\nonumber & = \bar{h}_1 \cdots \bar{h}_{j-1}z_1\bar{h}_{j+1}z_2\bar{h}_{j+2} \cdots \bar{h}_m
\end{align}
since $[\bar{h}_{j+1},z_2]$ is also an evaluation of $h_{j+1}$ and $z_1$ eradicates $f$ in $A$. Repeating this argument allows us to shift $z_2$ all the way to the right.  Consequently, $\bar{f}=0$, in this first case. \\

\noindent{Case (2):} Suppose that $x_i$ falls in $h_j=[x_{j_1},x_{j_2},\ldots,x_{j_k}]$, where $k\ge2$. Then part (1) of Lemma \ref{commutator} allows us to assume that $x_{j_1}=x_i$.
Thus, when we evaluate $x_i$ by $z_1z_2$, part (2) of Lemma 
\ref{commutator} yields
\begin{center}
$\bar{h}_j=[z_1z_2,a_{j_2},\ldots,a_{j_k}]=[z_1,z_2a_{j_2},\ldots,a_{j_k}]+[z_2,a_{j_2}z_1,a_{j_3}\ldots,a_{j_k}]$. 
\end{center}
It now follows that $\bar{f}=0$, as required.
\qed

\begin{defi}
We shall call a multilinear polynomial $f$ {\em distinctly proper} whenever it can be written in the form
\begin{center}
$f=g(h_1,\ldots,h_m)$,
\end{center}
where $g$ is multilinear and each polynomial $h_j$ is a simple Lie commutator of length at least $2$.
\end{defi}

\begin{thm}
\label{disproper}
Distinctly proper polynomials are marginally perfect.  In particular, derived series commutators are marginally perfect.
\end{thm}

\Proof Suppose $f$ is a distinctly proper polynomial.  Then 
\begin{center}
$f(x_1,\ldots,x_n)=g(h_1,\ldots,h_m)$,
\end{center}
where each $h_j$ is a simple Lie commutator of length at least $2$. Let $z_1,z_2$ be elements in $\widehat{\S}_A(f)$. 
We need to prove that $f$ evaluates to zero whenever some $x_i$ is evaluated to $z_1z_2$. Notice that $x_i$ must fall in exactly one 
$h_j=[x_{j_1},x_{j_2},\ldots,x_{j_k}]$, where $k\ge2$. Part (1) of Lemma \ref{commutator} and the multilinearity of $g$ allows us to assume that $h_j$ starts with the indeterminate $x_i$; that is, $x_{j_1}=x_i$.
Thus, when we evaluate $x_i$ by $z_1z_2$, part (2) of Lemma 
\ref{commutator} yields
\begin{center}
${\bar h}_j=[z_1z_2,a_{j_2},\ldots,a_{j_k}]=[z_1,z_2a_{j_2},\ldots,a_{j_k}]+[z_2,a_{j_2}z_1,a_{j_3}\ldots,a_{j_k}]$. 
\end{center}
It now follows from the multilinearity of $g$ that $f$ evaluates to zero, as required.
\qed

\begin{thm} If $f(x_1,\ldots,x_n)$ is distinctly proper, then $h=[f,x_{n+1}]$ is marginally perfect. In particular, all outer Lie commutators of the form $[g_n, x_{2^n+1}]$ are marginally perfect.
\end{thm}
\Proof
Let $A$ be an algebra, and that suppose that $z_1, z_2\in A$ eradicate $h$.  If we evaluate $x_i$ to $z_1z_2$, where $1\le i\le n$, then we obtain $\bar{h}=0$, exactly as in the proof of Theorem \ref{disproper}.  If we evaluate $x_{n+1}$ to $z_1z_2$, then 
$$\bar{h}=[\bar{f},z_1z_2]=z_1[\bar{f},z_2]+[\bar{f},z_1]z_2=0, $$
by part (1) of Lemma \ref{id}.
\qed

The following interesting special cases of Problem B remain open. {\em
\begin{itemize}
\item Are all outer Lie commutators marginally perfect or concise?
\item Are all Engel polynomials marginally perfect or concise?
\item Are all standard polynomials $s_n$ marginally perfect or concise?
\end{itemize}
}

\section{Corresponding verbal and marginal subspaces}

Hall raised the following question relating corresponding verbal and marginal subgroups:
{\em Let $\theta$ be a word in the free group such that the marginal subgroup, $\widehat{\theta}(G)$, of a group $G$ has finite index in $G$.  Does it follow that the corresponding verbal subgroup, $\theta(G)$, is finite?} He was inspired by a well-known result of Baer (\cite{RB}), which asserts that, for all positive integers $n$, if $G$ is a group such that $G/\centre_n(G)$ is finite, then the $(n+1)^\text{th}$ term of the lower central series of $G$ is also finite. Thus, Hall's question has positive solution when $\theta$ is the group commutator $(x_1,\dots,x_{n+1})$.  Hall proved that the converse of his question has a negative solution, in general, by constructing counterexamples to the converse of Baer's theorem.  See \cite{PH} for the original discussion of all these facts.

Stewart proved the following result in \cite{STW}, which can be considered as a positive solution to the algebra-theoretic analogue of Hall's problem. His proof works for all nonassociative algebras $A$, not just associative algebras.

\begin{thm} 
Let $A$ be an algebra, and suppose that $f$ is a polynomial with the property that $A/\widehat{\S}_A(f)$ is finite-dimensional. Then $\S_A(f)$ is also finite-dimensional. In other words, if $f$ is marginally a PI on $A$ then $f$ is almost a PI on $A$.\end{thm}

It is natural then to ask whether the converse of Theorem 4.1 holds: 

\begin{probC} Let $A$ be an algebra, and let $f$ be a polynomial in the free $K$-algebra.
If $\S_A(f)$ is finite-dimensional, under what conditions on $A$ and $f$ is $A/\widehat{\S}_A(f)$ necessarily finite-dimensional? Is this always true?
\end{probC}

For example, consider the case when $A$ is residually finite-dimensional and $f$ is any polynomial such that $\S_A(f)$ is finite-dimensional.  Then there exists an ideal $I$ of finite-codimension in $A$ such that $I\cap\S_A(f)=0$.  It follows that $I\subseteq \widehat{\S}_A(f)$, so that $A/\widehat{\S}_A(f)$ is finite-dimensional.  This proves:

\begin{prop} Let $f$ be a polynomial in the free $K$-algebra, and suppose that $A$ is a residually finite-dimensional algebra such that $\S_A(f)$ is finite-dimensional.  Then $A/\widehat{\S}_A(f)$ is finite-dimensional, too.
\end{prop}

However, the following construction taken from \cite{RS} shows that Problem C has a negative solution for general algebras $A$, even for `nice' polynomials $f$.  Recall that the lower central series of $A$, when viewed as a Lie algebra, is defined by $\g_1(A)=A$ and $\g_{n+1}(A)=[\g_{n}(A),A]$, for all $n\ge1$. 

\begin{eg} {\em (see Section 7 in \cite{RS})}
Let $K(\a)$ be any simple field extension of a base field $K$, and let $V$ be a vector space over $K(\a)$ with basis $\{v_1,\ldots,v_{n+1}\}$, for some fixed positive integer $n$. Now let $E$ denote the (non-unital) Grassmann-like $K(\a)$-algebra generated by $V$ subject to the relations
\begin{center}
$v_jv_i=\a v_iv_j$,
\end{center}
for all $1 \leq i \leq j \leq n+1$. Notice that these relations imply that $v_i^2=0$ except when $\a=1$; thus, in the case when $\a=1$, we impose the additional relations $v_i^2=0$, for each $i$. It is easy to see that $E$ has a $K(\a)$-basis consisting of all the monomials of the form
\begin{center}
$v_{i_1} \cdots v_{i_k}$,
\end{center}
where $1 \leq i_1 < \cdots < i_k \leq n+1$ and $1 \leq k \leq n+1$. Clearly $E^{n+2}=0$ and $\ann(E)=E^{n+1}=K(\a)v_1 \cdots v_{n+1}$. By induction, we also have
\begin{center}
$\ann^m(E)=E^{n+2-m}$ (see Example 3.2),
\end{center}
for each $0 \leq m \leq n+1$.

Now let $A$ be the algebra formed by identifying the elements corresponding to $v_1 \cdots v_{n+1}$ in each copy of the direct sum of countably many copies of $E$. Then the following statements hold.
\begin{enumerate}
\item If $\a=1$, then $A$ is a commutative $K$-algebra with the property that 
\begin{center}
$\dim_K(\S_A(x_1 \cdots x_{n+1}))=\dim_K(A^{n+1})=1$, 
\end{center}
and yet 
\begin{center}
$A/\widehat{\S}_A(x_1 \cdots x_{n+1})=A/\ann^n(A)=A/A^2$
\end{center}
is infinite-dimensional.
\item If $\a$ is a primitive root of unity whose order exceeds $n$, then $A$ has the properties that
\begin{align}
\nonumber \dim_K(\S_A([x_1,\ldots,x_{n+1}]))&=\dim_K(\g_{n+1}(A))\\
\nonumber &=\dim_K(K(\a))<\infty,
\end{align}
and yet
\begin{align}
\nonumber A/\widehat{\S}_A([x_1,\ldots,x_{n+1}])=A/\centre_n(A)=A/A^2,
\end{align}
is infinite-dimensional.
\end{enumerate}
\end{eg}

Notice that Example 4.3 shows that the converse of Theorem 4.1 fails for both associative and Lie algebras.  In closing, we remark that it remains conceivable that Problem C has a positive solution whenever $A$ is finitely generated.

\bigbreak \noindent {ERIC JESPERS}\\ { Department of
Mathematics}\\ {Vrije Universiteit Brussel} \\ {
Pleinlaan 2, 1050 Brussel}\\ {Belgium}\\ {e--mail:
\it Eric.Jespers@vub.be}

\bigbreak \noindent {DAVID RILEY}\\ { Department of
Mathematics}\\ {Western University}
\\ { London, Ontario, N6A 5B7}\\ {Canada}\\
{e--mail: \it DMRiley\@@uwo.ca}

\bigbreak \noindent {MAYADA SHAHADA}\\ { Department of
Mathematics and Statistics}\\ {Dalhousie University}
\\ { Halifax, Nova Scotia, B3H 4R2}\\ {Canada}\\
{e--mail: \it Mayada.Shahada\@@dal.ca}
\end{document}